\documentclass[11pt, reqno]{amsart}

 \usepackage{amsmath}
 \usepackage{amssymb}
\usepackage{bm}
\usepackage{color}

\newcommand{\mysection}[1]{\section{#1}
      \setcounter{equation}{0}}

\newtheorem{theorem}{Theorem}[section]
\newtheorem{lemma}[theorem]{Lemma}

\newtheorem{corollary}[theorem]{Corollary} 

\theoremstyle{definition}
\newtheorem{assumption}{Assumption}[section]

\theoremstyle{remark}
\newtheorem{remark}{Remark}[section]

\newcommand{\tr}{\text{\rm tr}\,}

\newcommand\bbeta{\text{\raise-.2ex\hbox{$\bm{\beta}$}}}
\newcommand\dist{{\rm dist}\,}

\makeatletter
 \def\dashint{%
 \operatorname%
 {\,\,\text{\bf--}\kern-.98em\DOTSI\intop\ilimits@\!\!}}
 \makeatother

\newcommand\bR{\mathbb{R}}

\newcommand\bS{\mathbb{S}}

\newcommand\bZ{\mathbb{Z}}

\newcommand\cL{\mathcal{L}}

\newcommand\frL{\mathfrak{L}}

\newcommand\supinf{\operatornamewithlimits{sup\,\,\,inf}}

\begin{document}

\title[On Isaacs equations]
{To the theory of viscosity solutions for
uniformly parabolic Isaacs equations}

\author{N.V. Krylov}
\thanks{The  author was partially supported by
 NSF Grant DMS-1160569}
\email{krylov@math.umn.edu}
\address{127 Vincent Hall, University of Minnesota,
 Minneapolis, MN, 55455}



\begin{abstract}
We show how a theorem about  the
solvability in $W^{1,2}_{\infty}$ of special
parabolic Isaacs equations
can be used to obtain the existence and uniqueness
of viscosity solutions of general
uniformly nondegenerate parabolic Isaacs equations.
We apply it also to establish the $C^{1+\chi}$
regularity of viscosity solutions and
show that finite-difference approximations
have an algebraic rate of convergence.
The main coefficients of the Isaacs equations
are supposed to be in $C^{\gamma}$ with
respect to the spatial variables with $\gamma$
slightly less than $1/2$.
 
\end{abstract}

\maketitle

\mysection{Introduction}

The goal of this article is to present
a purely PDE exposition of some key
results in the theory of viscosity
solutions for uniformly
nondegenerate parabolic Isaacs equations.
We did the same for elliptic equations
in \cite{Kr_15}, and this article is its natural
continuation.

Let $\bR^{d}=\{x=(x^{1},...,x^{d})\}$
be a $d$-dimensional Euclidean space.
Assume that we are given separable metric spaces
  $A$ and $B$,   and let,
for each $\alpha\in A$ and $\beta\in B$, 
  the following 
  functions on $\bR^{d+1}=\{(t,x):t\in\bR,x\in\bR^{d}\}$ be given: 

(i) $d\times d $
matrix-valued $a^{\alpha\beta} $,

(ii)
$\bR^{d}$-valued $b^{\alpha\beta} $, and

(iii)
real-valued  functions 
$c^{\alpha\beta}  \geq0$,   
  $f^{\alpha\beta}  $, and  
$g $. 

Let $\bS$ be the
set of symmetric $d\times d$ matrices, and for
$(u_{ij})\in\bS$, $(u_{i})\in\bR^{d}$, and $u\in\bR$
introduce 
$$
F(u_{ij},u_{i},u,t,x)=
\supinf_{\alpha\in A\,\,\beta\in B}
[a^{\alpha\beta}_{ij}( t,x)u_{ij}  +
b ^{\alpha\beta}_{i }(t, x)u_{i} -c^{\alpha\beta} (t, x)u
+f^{\alpha\beta}(t,x)],
$$
where and everywhere below the summation convention is enforced
and the summations are done inside the brackets.

For  sufficiently smooth functions $u=u(t,x)$ introduce
$$
L^{\alpha\beta} u(t,x)=a^{\alpha\beta}_{ij}( t,x)D_{ij}u(t,x)+
b ^{\alpha\beta}_{i }( t,x)D_{i}u(t,x)-c^{\alpha\beta} ( t,x)u(t,x),
$$
where, naturally, $D_{i}=\partial/\partial x^{i}$, $D_{ij}=D_{i}D_{j}$.
Denote
$$
F[u](t,x)=F(D^{2}u(t,x),Du(t,x),u(t,x),t,x)
$$
\begin{equation}
                                                     \label{1.16.1}
=\supinf_{\alpha\in A\,\,\beta\in B}
[L^{\alpha\beta} u(t,x)+f^{\alpha\beta} (t,x)],
\end{equation}
where $Du$ is the gradient of $u$ and $D^{2}u$ is its Hessian.

Also  
 fix an open bounded subset $G$ of $\bR^{d}$  with
$C^{2}$ boundary and $T\in(0,\infty)$. We denote 
the parabolic boundary of the cylinder
$Q=(0,T)\times  G $ by
$$
\partial'Q
=(\partial Q)\setminus (\{0\}\times G).
$$

Under appropriate assumptions,   which amount to the boundedness
and continuity with respect to $(t,x)$ of the 
coefficients and the free term and uniform
nondegeneracy of $a^{\alpha\beta}( t,x)$, the Isaacs equation
\begin{equation}
                                                     \label{4.1.4}
\partial_{t}u+F[u]=0,
\end{equation}
where $\partial_{t}=\partial/\partial t$,
in $Q$ with boundary condition $u=g$ on $\partial' Q$
has a   viscosity solution $v\in C(\bar{Q})$,
which is proved in Theorem 2.1 of 
Crandall,   Kocan,  Lions, and  \'Swi{\c e}ch
\cite{CKLS_99} (1999). Recall that $v\in C(\bar Q)$ is a viscosity
solution if for any $\phi(t,x)$ which is smooth  in $\bar Q$
and any point $(t_{0},x_{0})\in [0,T)\times G$
 at which $\phi-w$ attains

(i) 
a local (relative to $[0,T)\times G$)
 maximum which is zero it holds that $\partial_{t}\phi
(t_{0},x_{0})+F[\phi](t_{0}x_{0})\leq0$,

(ii) a local minimum
 which is zero we have $\partial_{t}\phi
(t_{0},x_{0})+F[\phi](t_{0}x_{0})\geq0$.

In \cite{CKLS_99} also the existence of minimal and maximal
continuous viscosity solution is proved. It turns out that
this fact holds true also for $L_{p}$-viscosity
solutions in the framework of the Isaacs
equations even if the coefficients and the free terms
are just measurable
(see, for instance, \cite{Kr13}).
Crandall,   Kocan,   and  \'Swi{\c e}ch
in \cite{CKS00} (2000) proved the existence of
continuous $L_{p}$-viscosity solutions
when the coefficients are only measurable with respect to $(t,x)$
and $\sup_{\alpha,\beta}|f^{\alpha\beta}|\in L_{p}$.
In the continuous case continuous $L_{p}$-viscosity solutions
are automatically viscosity solutions and this reproves
part of the results in \cite{CKLS_99}.

It seems that
in the parabolic case much less is known about uniqueness
 without convexity assumptions on $F$
with respect to $D^{2}u$ than in the elliptic case.
 One could derive uniqueness
for parabolic equations considering them as degenerate
elliptic ones and using the results
in  Crandall,   Ishii, and   Lions
  \cite{CIL_92} (1992), but this would require the coefficients
and the free term to be almost Lipschitz in $(t,x)$
(see Section 5.A there).
In Lemma 6.2 of \cite{CKS00} the
 uniqueness even of $L_{p}$-viscosity
solution is proved for equations which in our terms
have coefficients independent of $(t,x)$ and $f^{\alpha\beta}(t,x)
=f^{\alpha\beta} +f (t,x)$. We are going to prove uniqueness
under the assumption that $a^{\alpha\beta}$
are $\gamma$-H\"older continuous in $x$, with $\gamma<1/2$.
 The remaining
coefficients and the free term are assumed to be
uniformly continuous in $(t,x)$.
This result when  $a^{\alpha\beta}$ are Lipschitz
continuous only with respect to $x$ seems to be
available from \cite{Zh} (see Theorem 3.1 there).
In Corollary 3.5 of \cite{Zh} a comparison result,
and hence uniqueness,
is stated even for the case of just continuous 
$a^{\alpha\beta}$ (in our setting).

Concerning regularity of solutions note that the
 interior $C^{1+\chi}$-regularity
was established by Wang \cite{Wa92} (1992) under 
the assumption that 
$F$ is almost independent of $Du$ and (in our setting)
the coefficients are uniformly sufficiently  close
to the ones which are uniformly continuous with respect
to $(t,x)$ uniformly in $\alpha,\beta$.
Then Crandall,  Kocan, and \'Swi{\c e}ch \cite{CKS00}
(see there Theorem 7.3) generalized
the result of \cite{Wa92} to the case of full equation
 and continuous $L_{p}$-viscosity solutions assuming the
same kind of dependence on $(t,x)$ of $a^{\alpha\beta}(t,x)$
 as in \cite{Wa92}. As we know from \cite{Kr_14_0},
in the case of the Isaacs equations it is enough
to have $a^{\alpha\beta}(t,x)$ in VMO with respect to $x$.
We use this result in the form of Theorem \ref{theorem 3.31.1}.

 Last issue we are dealing with is the rate of convergence
of finite difference 
approximations of solutions of the Isaacs equations.
Caffarelli and Souganidis \cite{CS} (2008)
gave a method of establishing the rate of convergence
for elliptic equations of the form $F(D^{2}u(x))=f(x)$.
Turanova \cite{Tu_13_1} extended the results of 
\cite{CS} to $F$'s explicitly depending on $x$
and in \cite{Tu_14} considered parabolic equations
of the form $\partial_{t}u+F(D^{2}u)=0$.
We obtain the rate of convergence for the general 
parabolic Isaacs
equations under the same assumptions under which
we prove the existence and uniqueness of viscosity
solutions  
additionally assuming that the coefficients and the free
term are H\"older continuous with respect to
$(t,x)$ with a constant and exponent
independent of $\alpha,\beta$.

It is worth saying that in all the above cited papers
equations much more general than the Isaacs equations
are considered and we discuss their results in
the case of the Isaacs equations just to be able
to compare them with our results.
 Our methods are quite
different from the methods of the above cited
articles and are based on an approximation theorem
(Theorem \ref{theorem 3.30.1}, the proof of which, by the way,
does not even require any knowledge of the theory
of partial differential equations).

 The article is organized as follows.
In Section \ref{section 2.26.3} we present
our main results and prove all of them
apart from Theorem \ref{theorem 4.1.2},
 assertion (ii) of Theorem \ref{theorem 4.1.3},
and Theorem \ref{theorem 4.1.4},
which are proved in Sections \ref{section 4.3.3},
 \ref{section 4.3.4}, and \ref{section 5.3.1}, respectively, after
a rather long Section
\ref{section 4.3.2} containing a comparison
theorem for smooth functions.

Our equation are considered in $C^{2}$ cylinders
with $W^{1,2}_{\infty}$ boundary data. These
restrictions can be
considerably relaxed and we leave doing that to the interested
reader.
 
\mysection{Main results}
                                           \label{section 2.26.3}
 
Fix some constants $\delta\in(0,1)$ and $K_{0}\in[0,\infty)$.
 Set
$$
\bS_{\delta}=\{a\in\bS:\delta|\xi|^{2}\leq a_{ij}\xi_{i}\xi_{j}
\leq\delta^{-1}|\xi|^{2},\quad \forall\, \xi\in\bR^{d}\}. 
$$

In the following assumption the small parameter  $\chi\in(0,1)$,
 which depends only on
$\delta$ and $d$, is a constant  to be specified in
Theorem \ref{theorem 3.31.1} and
$$
\gamma=\frac{4-3\chi}{8-4\chi}\quad(<1/2).
$$

\begin{assumption}
                                    \label{assumption 1.9.1}
(i) The  functions $a^{\alpha\beta}$, $b^{\alpha\beta}$,
$c^{\alpha\beta}$, and $f^{\alpha\beta}$ 
are continuous with respect to
$\beta\in B$ for each $(\alpha, t,x)$ and continuous with respect
to $\alpha\in A$ uniformly with respect to $\beta\in B$
for each $(t,x)$, they are also uniformly continuous
with respect to $t$ uniformly with respect to $\alpha,\beta,x$,
the function $g$ is continuous and
$$
\|g\|_{ W^{1,2}_{\infty }(\bR^{d+1})}\leq K_{0},
$$

(ii) 
For all values of the arguments
$$
   |b^{\alpha\beta} |,
|c^{\alpha\beta} |,|f^{\alpha\beta} |
\leq K_{0},
$$

(iii) For any $(\alpha,\beta )\in A\times B $ and $(t,x),
(t,y)\in\bR^{d+1}$
we have
$$
\|a^{\alpha\beta}(t, x )-a^{\alpha\beta}( t,y )\|
\leq K_{0} |x-y|^{\gamma} ,
$$
$$
|u^{\alpha\beta}( t,x )-u^{\alpha\beta}(t, y )| 
\leq K_{0}\omega( |x-y|),
$$
where   $u=b,c,f$,  and $\omega$ is a fixed continuous
increasing 
function on $[0,\infty)$, such that $\omega(0)=0$,
and for a matrix $\sigma$ we denote $\|\sigma\|^{2}
=\tr\sigma\sigma^{*}$,

(iv) For all values of the
arguments $a^{\alpha\beta}\in\bS_{\delta}$.
 
\end{assumption}

\begin{remark}
The H\"older exponent 
 in Assumption \ref{assumption 1.9.1} (iii)
is certainly not sharp. For instance,
if $\chi$ is close to one the exponent $\gamma$ should
approach zero.

We wanted to reflect better what is going on when $\chi$
is close to zero, since nobody knows how much larger it is.
\end{remark}

In \cite{Kr13}
 a convex positive homogeneous of degree
one  function $P(u_{ij})$ is constructed
on $(u_{ij})\in\bS$
such that at all points of differentiability
of $P$ with respect to $(u_{ij})$ we have 
$(P_{u_{ij}})\in\bS_{\hat \delta}$,
where $\hat{\delta}$ is a constant in $(0,\delta)$
depending only on
$d $  and $\delta$.  This function
is constructed once only $d$ and $\delta$ are given
and possesses some additional properties to be mentioned and used
below. By $P[u]$ we denote $P(D_{ij}u)$.

The following result in what concern   equation
\eqref{3.30.6} is part of Theorem 2.1 of \cite{Kr13}.

\begin{theorem}
                                           \label{theorem 3.30.1}
For any $K\geq0$ each of the equations
\begin{equation}
                                                 \label{3.30.6}
\partial_{t}u_{ K}+\max(F[u_{  K}],  P[  u
_{ K}]- K)=0,
\end{equation}
\begin{equation}
                                                 \label{3.30.06}
\partial_{t}u_{-K}+\min(F[u_{-K}],-P[- u
_{- K}]+K)=0,
\end{equation}
in $Q$ with boundary condition $u_{\pm K}=g$ on $\partial' Q$
has a unique solution in the class $W^{1,2}_{\infty,loc}(Q)
\cap C(\bar Q)$. Furthermore, there is a constant $N$
such that for any $K\geq0$, we have 
$|\partial_{t}u_{\pm K}|$,
$\rho|D^{2}u_{\pm K}|$, $|Du_{\pm K}|\leq N(1+K)$ in $Q$
 (a.e.), where
$\rho(t,x)=\dist(x,G^{c})$.

\end{theorem}

The result concerning equation \eqref{3.30.06}
 is obtained just by observing that
the operator $-F[-u]$ fits into the scheme of \cite{Kr13}
where equations more general than the Isaacs equations are
treated.

The functions $u_{\pm K}$ are the central objects of our
investigation. Here is a simple property they
possess.

\begin{lemma}
                                 \label{lemma 4.1.1}
There exists a constant $N$, depending only on
$d,\delta$,  $K_{0}$, and $G$,  such that in $Q$
$$
|u_{\pm K}-g| \leq N\rho,\quad
|u_{\pm K}| \leq N.
$$

\end{lemma}

This result for $u_{K}$
follows from the fact that 
  $|\max(F[0],-K)|\leq|F[0]|$, $g\in 
W^{1,2}_{\infty }(\bR^{d+1})$, $G\in C^{2}$, and
$u_{K}$ satisfies a linear equation
$$
\partial_{t}u_{K}+
a_{ij}D_{ij}u_{K}+b_{i}D_{i}u_{K}-cu_{K}+f=0,
$$
where $(a_{ij})\in\bS_{\hat \delta}$, $|(b_{i})|\leq K_{0}$,
$K_{0}\geq c\geq 0$, and $|f|\leq K_{0}$.
The case of $u_{-K}$ is quite similar.

To characterize some smoothness properties of $u_{\pm K}$
introduce $C^{1+\chi}(Q)$ as the space 
of functions on $Q$ continuously differentiable
with respect to $x$ and
with finite norm given by
$$
\|u\|_{C^{1+\chi}(Q)}=\sup_{Q}|u|
+\sup_{(t,x),(t,y)\in Q}\frac{|u(t,x)-u(t,y)|}{|x-y|}
+[u]_{C^{1+\chi}(Q)},
$$
where
$$
 [u]_{C^{1+\chi}(G)}=
 \sup_{(t,x),(s,x)\in Q}\big[\frac{|u(t,x)-u(s,x)|}
{|t-s|^{(1+\chi)/2}}
$$
$$
+ \frac{|Du(t,x)-Du(t,y)|}{|x-y|^{\chi}}
+ \frac{|Du(t,x)-Du(s,x)|}{|t-s|^{\chi/2}}\big]
$$
(the last term in the brackets can be dropped
yielding an equivalent norm
as follows from simple interpolation inequalities).

For $\varepsilon>0$ introduce
$$
  G_{\varepsilon}=\{x\in G:
\dist (x,\partial G)>\varepsilon\},
\quad Q_{\varepsilon}=(0,T-\varepsilon^{2})
\times G_{\varepsilon}.
$$

The following result in what concern $u_{K}$ is a consequence
of Theorem 5.4 of \cite{Kr_14_0}. The case of $u_{-K}$
is treated as in Theorem \ref{theorem 3.30.1}.

\begin{theorem}
                                            \label{theorem 3.31.1}
There exists a constant $\chi\in(0,1)$,
depending only on
$\delta$ and $d$, and there exists a constant 
$N$, depending only on $K_{0}$,
$\delta$, $d$, and $G$ such that for any $\varepsilon>0$
(such that $G_{\varepsilon}\ne\emptyset$)
\begin{equation}
                                                 \label{3.30.8}
 \|u_{\pm K}\|_{C^{1+\chi}(Q_{\varepsilon})}\leq
N\varepsilon^{-1-\chi}.
\end{equation}

\end{theorem}

Estimate \eqref{3.30.8} without specified dependence
on $\varepsilon$ of the right-hand side follows from
  Theorem 7.3 of \cite{CKS00}. We use \eqref{3.30.8}
as it is stated in the proof of 
the following result which is central in the paper.
Fix a constant $\tau\in(0,1)$.

\begin{theorem}
                                            \label{theorem 4.1.2}
For $K\to\infty$ we have $|u_{K}-u_{-K}|\to0$
uniformly in $Q$. Moreover, if
\begin{equation}
                                                    \label{4.5.1}
\omega(t)=t^{\tau},
\end{equation}
then there exist    constants $\xi\in(0,1)$,
depending only on $\tau$,
    $d$, $K_{0}$, and $\delta$,
and   $N\in(0,\infty)$,   depending only on  
  $\tau$,  $d$, $K_{0}$,  $\delta$, and   
$Q$,
such that, if $K\geq 1 $, then
\begin{equation}
                                    \label{3.30.01}
|u_{K}-u_{-K}|\leq  NK^{-\xi}
\end{equation}
 in $Q $.

\end{theorem}
We prove Theorem \ref{theorem 4.1.2}
 in Section \ref{section 4.3.3} by adapting 
an argument
from Section 5.A of \cite{CIL_92} explaining how to
prove the comparison principle for
$C^{1+\chi}$ subsolutions and supersolutions.
We use a quantitative and parabolic version of this argument.

\begin{theorem}
                                            \label{theorem 4.1.3}
(i) The limit
$$
v:=\lim_{K\to\infty}u_{K}
$$
exists,

(ii) The function $v$ is a unique continuous in $\bar Q$
viscosity solution
of \eqref{4.1.4} with boundary condition $v=g$ on
$\partial' Q$,

(iii) If condition \eqref{4.5.1} is satisfied, then
for   $K\geq1$ we have $|u_{K}-v|\leq  NK^{-\xi}$, 

(iv) For any $\varepsilon\in(0,1]$
(such that $G_{\varepsilon}\ne\emptyset$)
$$
 \|v\|_{C^{1+\chi}(Q_{\varepsilon})}\leq
N\varepsilon^{-1-\chi},
$$
where $N$ is the constant from \eqref{3.30.8}.

\end{theorem}

Assertions (i), (iii), and (iv) are simple consequences
of Theorems \ref{theorem 3.31.1}
and \ref{theorem 4.1.2} and the maximum principle.
Indeed,
notice that $\partial_{t}u_{K}+F[u_{K}]\leq0$ and $
\partial_{t}u_{-K}+F[u_{-K}]\geq0$.
Hence by the maximum principle $u_{K}\geq u_{-K}$.
Furthermore, again by the maximum principle
$u_{K}$ decreases and $u_{-K}$ increases as $K$
increases.  This takes care
of assertions (i), (iii), and (iv).

Assertion (ii) is proved in Section \ref{section 4.3.4}.

To state one more result we introduce some new objects.
As is well known (see, for instance, \cite{KT92}),
there exists a finite set $\Lambda=\{l_{1},...,l_{d_{2}}\}
\subset\bZ^{d}$  containing all vectors from the standard 
orthonormal basis of $\bR^{d}$
such that one has the following representation
$$
L^{\alpha\beta}u(t,x)=
a^{\alpha\beta}_{k}(t,x)D_{l_{k}}^{2}u(t,x)+
\bar b^{\alpha\beta}_{k}(t,x)D_{l_{k}}u(t,x)
-c^{\alpha\beta}(t,x)u(t,x),
$$
where $D_{l_{k}}u(x)=\langle D u ,l_{k}\rangle$,
$a^{\alpha\beta}_{k}$ and $\bar b^{\alpha\beta}_{k}$ are certain
bounded functions and $a^{\alpha\beta}_{k}\geq\delta_{1}$,
with a constant $\delta_{1}>0$. One can even arrange for
such representation to have the coefficients 
$a^{\alpha\beta}_{k}$ and $\bar b^{\alpha\beta}_{k}$ 
with the same regularity
properties with respect to $(t,x)$ as the original ones
$a^{\alpha\beta}_{ij}$ and $b_{i}^{\alpha\beta}$
(see, for instance, Theorem 3.1 in \cite{Kr_11}).
Define
$B$ as the smallest closed ball centered
at the origin containing $\Lambda$,
and for $h>0$ set $\bZ^{d}_{h}=h\bZ^{d}$,
$$
G_{(h)}=G\cap\bZ^{d}_{h},\quad Q_{(h)}=((0,T)\cap \bZ_{h^{2}})
\times G_{(h)},
$$
$$
Q^{o}_{(h)}=\{(t,x)\in Q_{ 
(h)}:x+hB\subset G, t+h^{2}<T\},\quad
\partial'_{h}Q=Q_{(h)}\setminus Q^{o}_{(h)}.
$$

Next, for $h>0$ we introduce
$$
\delta_{h,t}u(t,x)=\frac{u(t+h^{2},x)-u(t,x)}{h^{2}},\quad
\delta_{h,l_{k}}u(t,x)=\frac{u(t,x+hl_{k})-u(t,x)}{h} ,
$$
$$
\Delta_{h,l_{k}}u(t,x)=\frac{u(t,x+hl_{k})-
2u(t,x)+u(t,x-hl_{k})}{h^{2}} ,
$$
$$
L^{\alpha\beta}_{h}u(t,x)=
a^{\alpha\beta}_{k}(t,x)\Delta_{h,l_{k}}u(t,x)+
\bar b^{\alpha\beta}_{k}(t,x)\delta_{h,l_{k}}u(t,x)
-c^{\alpha\beta}(t,x)u(t,x),
$$
$$
F_{h}[u](t,x)=
\supinf_{\alpha\in A\,\,\beta\in B}[L^{\alpha\beta}_{h}u(t,x)
+f^{\alpha\beta}(t,x)].
$$

It is a simple fact, which may be shown as in \cite{KT92},
that   for each sufficiently small $h$ there exists a unique
  function $v_{h}$ on $Q_{(h)}$ such that
$$
\delta_{h,t}v_{h}+F_{h}[v_{h}]=0
$$
 on $Q^{o}_{(h)}$ and 
$v_{h}=0$ on $\partial'_{h}Q$.

Here is the result we were talking about above
and which is proved in Section \ref{section 5.3.1}.

\begin{theorem}
                                          \label{theorem 4.1.4}
 
Suppose that there exists $\gamma_{t}\in(0,1)$
such that for any $(\alpha,\beta )\in A\times B $ and $(t,x),
(s,x)\in\bR^{d+1}$
we have
$$
|u^{\alpha\beta}( t,x )-u^{\alpha\beta}(s, x )| 
\leq K_{0} |t-s|^{\gamma_{t}},
$$
where   $u=a,b,c,f$.
Let condition \eqref{4.5.1} be satisfied
and let $g=0$.
Then there exist constants $N$ and $\eta>0$ such that
for all sufficiently small $h>0$ we have on $Q_{(h)}$ that
$$
|v_{h}- v|\leq Nh^{\eta}.
$$
\end{theorem}

In this theorem we deal with implicit finite-difference
schemes and with the time step size rigidly related
to the space step size. This is actually irrelevant
and explicit or mixed schemes and different step sizes
can be treated in the same way. We do not do this just
not to overburden the exposition with lots of details of
relatively minor importance. The same should be said
about zero boundary data.

\mysection{An auxiliary result}
                                           \label{section 4.3.2}

In the following theorem   
 $G$ can be just any bounded domain.
Below by $C^{1,2}( Q)$ we mean the space of
functions $u=u(t,x)$ which are bounded and continuous
in $[0,T)\times G$ along with their derivatives $\partial_{t}u$,
$D_{ij}u$, $D_{i}u$.

\begin{theorem}
                                               \label{theorem 3.27.1}
Let  
 and $u,v\in C^{1,2}( Q)\cap C(\bar Q)$
be such that for a constant $K\geq1$
\begin{equation}
                                                 \label{3.30.5}
\partial_{t}u+\max(F[u],P[u]-K)\geq0\geq
\partial_{t}v+\min(F[v],-P[-v]+K)
\end{equation}
in $Q$ and $v\geq u$
on $\partial' Q$. Also assume that, for a constant $M \in[1,\infty)$, 
\begin{equation}
                                                 \label{3.28.6}
\|u,v\|_{C^{1+\chi}(Q)}\leq M .
\end{equation}
Then there exist  a  constant    $N \in(0,\infty)$,
depending only on $\tau$, the diameter of $G$,
    $d$, $K_{0}$, and $\delta$,
 and a  constant  $\eta\in(0,1)$,
depending only on $\tau$,
    $d$,   and $\delta$,
such that, if $K\geq NM^{1/\eta}$
and
\begin{equation}
                                     \label{5.27.1}
K\geq T^{-1},\quad G_{2/\sqrt{K}}\ne\emptyset,
\end{equation}  
then  in $Q$
\begin{equation}
                                                 \label{3.30.1}
u(t,x)-v(t,x)\leq  NK^{-\chi/4}
+ NM\omega(M^{-1/\tau}K^{-1}).
\end{equation}

\end{theorem}

\begin{remark}
                                            \label{remark 4.5.2}
Observe that for $\omega=t^{\tau}$ estimate
\eqref{3.30.1} becomes $u-v\leq  NK^{-\chi/4}+NK^{-\tau}$.
\end{remark}

To prove this theorem,
we are going to adapt to our situation an argument
from Section 5.A of \cite{CIL_92}. For that
we need a construction and two lemmas.
 From the start 
throughout the section we will only concentrate
on $K$ satisfying \eqref{5.27.1}.

First we introduce $\psi\in C^{2}(\bR^{d})$ as a global barrier for $G$,
that is, in $G$ we have $\psi\geq1$ and
$$
a_{ij}D_{ij}\psi+b_{i}D_{i}\leq-1
$$
for any $(a_{ij})\in \bS_{\hat \delta}, |(b_{i})|\leq K_{0}$.
Such a $\psi$ can be found in the form
$  \cosh \mu R-\cosh\mu|  x|$ for sufficiently large $\mu$
and $R$.

Then we take and fix a radially symmetric
with respect to $x$ function $\zeta=\zeta(t,x)$ of class
$ C^{\infty}_{0}(\bR^{d+1})$ with support in $(-1,0)\times B_{1}$.
For $\varepsilon>0$ we define $\zeta_{\varepsilon}(t,x)
=\varepsilon^{-d-2}\zeta(\varepsilon^{-2}t,\varepsilon^{-1}x)$
and for locally summable $u(t,x)$ introduce
\begin{equation}
                                                   \label{5.3.2}
u^{(\varepsilon)}(t,x)=u(t,x)*\zeta_{\varepsilon}(t,x).
\end{equation}
Recall some standard properties of parabolic mollifiers
in which no regularity properties of $G$ are required:
 If $u\in C^{ 1+\chi}(Q)$, then in $Q_{\varepsilon}$
$$
\varepsilon^{-1-\chi}|u-u^{(\varepsilon)}|
+\varepsilon^{ -\chi}|Du-Du^{(\varepsilon)}|
\leq N\|u
\|_{C^{ 1+\chi}(Q)},
$$
$$
 |u^{(\varepsilon)}
 |+|Du^{(\varepsilon)}| +\varepsilon^{1-\chi}|D^{2}u^{(\varepsilon)}|
 +\varepsilon^{1-\chi}|\partial_{t}u^{(\varepsilon)}|
$$
\begin{equation}
                                                    \label{3.31.1}
 +\varepsilon^{2-\chi}|D^{3}u^{(\varepsilon)}|
+\varepsilon^{2-\chi}|D\partial_{t}u^{(\varepsilon)}|
+\varepsilon^{3-\chi}|\partial_{t}D^{2}u^{(\varepsilon)}|
\leq N\|u
\|_{C^{ 1+\chi}(Q)}.
\end{equation}

Define the functions
$$
\bar u=u/\psi,\quad\bar v=v/\psi. 
$$
Replacing $M $ with $NM$, if necessary, where $N$
depends only on
 $d$ and the diameter of $G$, we may assume that
\eqref{3.28.6} holds with $\bar u,\bar v$ in place of $u,v$.

Next take  constants $\nu,\varepsilon_{0}\in(0,1)$,
introduce 
$$
\varepsilon=\varepsilon_{0} K^{-(1-\gamma)/(2\gamma)},
$$
and
consider the function 
$$
W(t,x, y)=\bar u(t,x)-\bar u^{(\varepsilon)}(t,x)-
[\bar v(t,y)-\bar u^{(\varepsilon)}(t,y)]-\nu K
 |x-y|^{2} 
$$
for $ (t,x),(t,y)\in \bar Q_{\varepsilon}$.
  Note that $Q_{\varepsilon}\ne
\emptyset$ and even $G_{2\varepsilon}\ne\emptyset$
owing to \eqref{5.27.1} and the fact that
$1-\gamma>\gamma$.

Denote by  $(\bar t, \bar x,  \bar y)$ a maximum point
of $W$ in $ [0,T-
\varepsilon^{2}] \times \bar G _{\varepsilon}^{2}$.
Observe that, obviously,
$$
\bar u(\bar t,\bar x)-\bar u^{(\varepsilon)}(\bar t,\bar x)-
[\bar v(\bar t,\bar y)-\bar u^{(\varepsilon)}(\bar t,\bar y)]-\nu K
 |\bar x-\bar y|^{2} 
$$
$$
\geq
\bar u(\bar t,\bar x)-\bar u^{(\varepsilon)}(\bar t,\bar x)-
[\bar v(\bar t,\bar x)-\bar u^{(\varepsilon)}(\bar t,\bar x)],
$$
which implies that
\begin{equation}
                                                    \label{4.21.2}
 |\bar x-\bar y|\leq N M/(\nu K) .
\end{equation}
where and below by $N$  with indices or without them we denote 
various constants
depending only on 
    $d$, $K_{0}$, $\delta$, and  the diameter of $G$,
unless specifically stated otherwise.
By the way recall that $\chi$ and, hence, $\gamma$
 depend  only on $d$ and $\delta$.

\begin{lemma}
                                              \label{lemma 3.30.1}
There exist  a constant   $\nu\in(0,1)$,
depending only on
    $d$, $K_{0}$,  $\delta$, and the diameter of $G$, and 
a constant $N\in(0,\infty)$
  such that, if
\begin{equation}
                                            \label{3.31.09}
K\geq N\varepsilon_{0}^{(\chi-1)/\eta_{1}}M^{1/\eta_{1}},
\end{equation}
where $  \eta_{1}=1-(1-\chi)(1-\gamma)/(2\gamma)$ ($>0$),
and $ \bar x ,\bar y \in G _{\varepsilon}$ 
and $ \bar t<T-\varepsilon^{2}$, then

(i) We have
\begin{equation}
                                            \label{3.27.6}
2\nu K|\bar x-\bar y|\leq NM \varepsilon^{\chi},
\quad|\bar x-\bar y|\leq\varepsilon/2,
\end{equation}

(ii) For any $\xi,\eta\in\bR^{d}$
\begin{equation}
                                            \label{3.27.7}
D_{ij}[\bar u -\bar u^{(\varepsilon)} ](\bar t,\bar x)\xi^{i}\xi^{j}
-D_{ij}[\bar v -\bar u^{(\varepsilon)} ](\bar t,\bar y)\eta^{i}\eta^{j}
\leq 2\nu K|\xi-\eta|^{2},
\end{equation}
\begin{equation}
                                            \label{3.31.2}
\partial_{t}[\bar u -\bar u^{(\varepsilon)} ]
(\bar t,\bar x) \leq 
 \partial_{t}[\bar v -\bar u^{(\varepsilon)} ]
(\bar t,\bar y) ,
\end{equation}

(iii) We have
\begin{equation}
                                               \label{3.27.8}
\partial_{t}\bar u(\bar t,\bar x)+
\supinf_{\alpha\in A\,\,\beta\in B}\big[a_{ij}^{\alpha\beta}
D_{ij}\bar u+\hat b^{\alpha\beta}_{i}D_{i}\bar u-
\hat c^{\alpha\beta}\bar u+
\hat f^{\alpha\beta}\big](\bar t,\bar x)\geq0,
\end{equation}

\begin{equation}
                                               \label{3.28.1}
\partial_{t}\bar v(\bar t,\bar y)+
\supinf_{\alpha\in A\,\,\beta\in B}\big[a_{ij}^{\alpha\beta}
D_{ij}\bar v+\hat b^{\alpha\beta}_{i}D_{i}\bar v-
\hat c^{\alpha\beta}\bar v+\hat f^{\alpha\beta}\big]
(\bar t,\bar y)\leq0,
\end{equation}
where
$$
 \hat b^{\alpha\beta}_{i}=  b^{\alpha\beta}_{i} 
+2a_{ij}^{\alpha\beta}\psi ^{-1}D_{j}\psi,\quad
\hat c^{\alpha\beta}=-\psi ^{-1}L^{\alpha\beta}\psi,
\quad \hat f^{\alpha\beta}=\psi ^{-1}f^{\alpha\beta}.
$$

\end{lemma}

Proof. The first inequality in
\eqref{3.27.6}  follows from \eqref{3.31.1} and
the fact that the first derivatives
of $W$ with respect to $x$ vanish  at $\bar x$, that is,
$D(\bar u-\bar u^{(\varepsilon)}
)(\bar x)=2\nu K(\bar x-\bar y)$.
Also the matrix of second-order derivatives of $W$
with respect to $(x,y)$
is nonpositive at $( \bar t,\bar x,\bar y)$ as well as its 
(at least one sided if $\bar t=0$)
derivative
with respect to   $t$, which yields (ii).

By taking $\eta=0$ in \eqref{3.27.7} and using the fact that
$|D^{2}\bar u^{(\varepsilon)}|\leq N M
\varepsilon^{\chi-1}$ we see that
$$
D^{2}\bar u(\bar t,\bar x)\leq 2\nu K+NM
\varepsilon^{\chi-1}.
$$
Furthermore 
$$
D_{ij} u=\psi D_{ij} \bar u+(D_{i}\psi)D_{j}\bar u
+(D_{i}\bar u)D_{j}\psi+(D_{ij}\psi)\bar u,
$$
which implies that
$$
D^{2}  u(\bar t,\bar x)\leq N(\nu K+ M
\varepsilon^{\chi-1}) .
$$
Similarly,
$$
D^{2}  v(\bar t,\bar x)\geq -N(\nu K+ M
\varepsilon^{\chi-1}) ,
$$
which yields
$$
P[u](\bar t,\bar x)
\leq N_{1}(\nu K+ M
\varepsilon^{\chi-1} ),\quad -P[-v](\bar t,\bar y)
\geq -
 N_{1}(\nu K+ M
\varepsilon^{\chi-1} ),
$$
\begin{equation}
                                                 \label{5.19.1}
F[u](\bar t,\bar x)
\leq N_{1}(\nu K+ M
\varepsilon^{\chi-1} ),\quad F[v](\bar t,\bar y)
\geq -
 N_{1}(\nu K+ M
\varepsilon^{\chi-1} ).
\end{equation}
Also it follows from \eqref{3.31.2} and \eqref{3.31.1} that
\begin{equation}
                                                 \label{5.19.2}
\partial_{t}u(\bar t,\bar x)\leq \partial_{t}v(\bar t,\bar y)
+N_{2} M
\varepsilon^{\chi-1} .
\end{equation}

Now, if $F[u](\bar t,\bar x)\leq P[u](\bar t,\bar x)-K$,
then at $(\bar t,\bar x)$
$$
0\leq
\partial_{t}u 
+\max(F[u] , P[u] -K)\leq\partial_{t}u+
N_{1}(\nu K+ M
\varepsilon^{\chi-1} )-K,
$$
$$
\partial_{t}u \geq K-N_{1}(\nu K+ M
\varepsilon^{\chi-1} )
$$
and at $(\bar t,\bar y)$
$$
0\geq\partial_{t}v 
+\min(F[v] , -P[-v]+K)
\geq K-2N_{1}(\nu K+ M
\varepsilon^{\chi-1} )-N_{2} M
\varepsilon^{\chi-1}.
$$
 This is impossible if we
  choose and fix $\nu$  such that
\begin{equation}
                                            \label{3.28.2}
(4N_{1} +N_{2})\nu\leq1/2,
\end{equation}
since, as
 is easy to see, $M\varepsilon^{\chi-1}\leq \nu K$ for $K$
satisfying \eqref{3.31.09}  and appropriate 
$N $, so that
$$
2N_{1}(\nu K+ M
\varepsilon^{\chi-1} )+N_{2} M
\varepsilon^{\chi-1}\leq 4N_{1}\nu K+N_{2}\nu K\leq K/2.
$$
It follows that $F[u](\bar t,\bar x)\geq P[u](\bar t,\bar x)-K$,
$$
\partial_{t}u(\bar t,\bar x)+F[u](\bar t,\bar x)\geq0,
$$
and this proves \eqref{3.27.8}.

Similarly, if $-P[-v](\bar t,\bar y)+K\leq F[v](\bar t,\bar y)$,
then at $(\bar t,\bar y)$
$$
0\geq\partial_{t} v -
 N_{1}(\nu K+ M
\varepsilon^{\chi-1} )+K,
$$
and at $(\bar t,\bar x)$ we have $\partial_{t} u\leq
N_{1}(\nu K+ M
\varepsilon^{\chi-1} ) -K+N_{2} M
\varepsilon^{\chi-1}$ and
$$
0\leq 
\partial_{t} u+\max(F[u],P[u]-K)\leq
-K+ 2N_{1}(\nu K+ M
\varepsilon^{\chi-1} )+N_{2} M
\varepsilon^{\chi-1},
$$
which again is impossible with the above choice of $\nu$
for $K$ satisfying \eqref{3.31.09}. This yields
\eqref{3.28.1}.

Moreover, not only $M\varepsilon^{\chi-1}\leq \nu K$ for
$K$
satisfying \eqref{3.31.09}, but we also have
$NM\varepsilon^{\chi-1}\leq \nu K$, where $N$ is taken from
\eqref{3.27.6}, if we increase $N$ in \eqref{3.31.09}.
This yields the second inequality in \eqref{3.27.6}.
 
 The lemma is proved.

\begin{lemma}
                                            \label{lemma 4.23.1}
For any $\mu\geq0$,
there exist   a  constant  $\eta>0$,
depending only on  
    $d$  and $\delta$, and a constant $N$ 
such that, if $K\geq NM^{1/\eta}$ and
\begin{equation}
                                                 \label{3.28.5}
W(\bar t,\bar x, \bar y)\geq 2K^{-\chi/4}
+\mu M\omega(M^{-1/\tau}K^{-1}),
\end{equation}
then 
\begin{equation}
                                            \label{3.27.4}
\bar u(\bar t,\bar x)- 
 \bar v(\bar t,\bar y) -
\nu K|\bar x-\bar y|^{2} 
\geq K^{-\chi/4}+\mu M\omega(M^{-1/\tau}K^{-1}).
\end{equation}
Furthermore,
$ \bar x,\bar y \in G_{2\varepsilon}$ and
$\bar t <T-\varepsilon^{2}$.
\end{lemma}

Proof.  
It follows from \eqref{4.21.2} 
that  (recall that $\nu$ is already fixed)
 \begin{equation}
                                                 \label{4.21.5}
|\bar u^{(\varepsilon)}(\bar t,\bar x)-
\bar u^{(\varepsilon)}(\bar t,\bar y)|\leq M  |\bar x-\bar y|
\leq N M^{2}/  K.
\end{equation}
Hence we have from \eqref{3.28.5} that  
$$
\bar u(\bar t,\bar x)- 
 \bar v(\bar t,\bar y) -
\nu K|\bar x-\bar y|^{2}\geq 2 K^{-\chi/4}-N M^{2}/  K
+\mu M\omega(M^{-1/\tau}K^{-1}),
$$
and \eqref{3.27.4} follows provided that 
\begin{equation}
                                                   \label{4.22.2}
NM^{2}/K\leq (1/4) K^{-\chi/4},
\end{equation}
which indeed holds if
\begin{equation}
                                                   \label{3.31.6}
K\geq N 
M^{1/\eta_{2}},
\end{equation}
 with $\eta_{2}=(4-\chi)/8>0$.
Here if $\bar x$ or $\bar y$ are 
 in $\bar G_{\varepsilon}
\setminus G_{2\varepsilon}$, then
for appropriate $\hat x\in\partial G $ and $\hat y\in\partial G $
either 
$$
\bar u(\bar t,\bar x)- 
 \bar v(\bar t,\bar y)\leq 2M \varepsilon+ 
\bar v(\bar t,\hat x)-\bar v(\bar t,\bar y)
\leq M (4\varepsilon+|\bar x-\bar y|)
$$
or
$$
\bar u(\bar t,\bar x)- 
 \bar v(\bar t,\bar y)\leq \bar u(\bar t,\bar x)
-\bar u(\bar t,\hat y)
+ \bar v(\bar t,\hat y)-\bar v(\bar t,\bar y)
\leq M (4\varepsilon+|\bar x-\bar y| ).
$$

In any case in light of \eqref{3.27.4}
\begin{equation}
                                                   \label{4.23.5}
4\varepsilon M +N M^{2}/K-
\nu K|\bar x-\bar y|^{2} \geq K^{-\chi/4}.
\end{equation}
Notice that (we need the following with a constant $N$
for the future)
\begin{equation}
                                                   \label{4.22.1}
NM \varepsilon\leq  (1/4) K^{-\chi/4}
\end{equation}
 for 
\begin{equation}
                                                   \label{3.31.7}
K\geq (4NM)^{1/\eta_{3}},
\end{equation}
where 
 $\eta_{3}:=(1-\gamma)/(2\gamma)-\chi/4 >0$. This and
\eqref{4.22.2} show that \eqref{4.23.5} is 
impossible   for 
\begin{equation}
                                                     \label{3.31.8}
 K\geq N (M^{1/\eta_{2}}
+M^{1/\eta_{3}}).
\end{equation}
Below we assume \eqref{3.31.8} and   for such $K$,
we have  $ \bar x,\bar y \in G_{\varepsilon}$.

Furthermore, if $\bar t= T-\varepsilon^{2}$, then
(recall \eqref{4.21.5} and that $\bar u\leq\bar v$
on $\partial' Q$)
$$
W(\bar t,\bar x, \bar y)\leq 
NM\varepsilon^{1+\chi}+\bar v(\bar t,\bar x)-
\bar v(\bar t,\bar y)+NM^{2}/K
$$
$$
\leq NM\varepsilon+M |\bar x-\bar y| 
 +NM^{2}/K\leq NM^{2}/K+NM\varepsilon ,
$$
which is less than $K^{-\chi/4}$ (cf.~\eqref{4.22.1} 
and \eqref{4.22.2}).
This is impossible due to \eqref{3.28.5}.
Hence, $\bar t <T-\varepsilon^{2}$ and this
finishes the proof of the lemma.

{\bf Proof of Theorem \ref{theorem 3.27.1}}.
Fix a (large) constant $\mu>0$ to be specified later
as a constant, depending only on  $d$, $K_{0}$,  $\delta$,
and the diameter of $G$,
take $\nu$ from Lemma \ref{lemma 3.30.1} and
first assume that 
$$
W(t,x ,y)\leq 2K^{-\chi/4}+\mu M\omega(M^{-1/\tau}K^{-1})
$$
for  $ (t,x),(t,y)\in \bar Q_{\varepsilon}$. Observe that for 
any point $(t,x)\in \bar  Q $
one can find a point $(s,y)\in \bar Q _{\varepsilon}$
with $|x-y|\leq\varepsilon$ and $|t-s|\leq\varepsilon^{2}$ and
then 
$$
\bar u(t,x)-\bar v(t,x)\leq \bar u(s,y)-\bar v(s,y)
+NM \varepsilon
$$
$$
\leq W(s,y, y)+NM \varepsilon
\leq 2K^{-\chi/4}+NM \varepsilon
+\mu M\omega(M^{-1/\tau}K^{-1}).
$$
In that case, owing to \eqref{4.22.1}, \eqref{3.30.1} holds for 
$K$ satisfying \eqref{3.31.7}.

It is clear now that, to prove the theorem, it suffices to find  
  $\mu$ such that the inequality  \eqref{3.28.5}
is impossible
if $K\geq NM^{1/\eta}$ with $N$
and $\eta$ as in the statement of the theorem.
 Of course, we will argue by 
contradiction and suppose that \eqref{3.28.5} holds.

Upon combining   Lemmas \ref{lemma 3.30.1} and 
\ref{lemma 4.23.1} 
we see that there exist  $\eta_{0}\in(0,1)$,
depending only on  
    $d$  and $\delta$, and $N\in(0,\infty)$
 such that, if
\begin{equation}
                                                     \label{3.30.3}
K\geq N\varepsilon_{0}^{(\chi-1)/\eta_{0}}M^{1/\eta_{0}},
\end{equation}
then 
\eqref{3.27.4} holds,
   $ \bar x ,\bar y \in G _{\varepsilon}$ 
and $\bar t<T-\varepsilon^{2}$, so that
we can use   the conclusions of Lemma \ref{lemma 3.30.1}.

By   denoting $\sigma^{\alpha\beta}=(a^{\alpha\beta})^{1/2}$
we may write
$$
a_{ij}^{\alpha\beta} =\sigma^{\alpha\beta}_{ik}
\sigma^{\alpha\beta}_{jk},
$$
and then \eqref{3.27.7} for $\xi^{i}=
\sigma^{\alpha\beta}_{ik}(\bar t,\bar x)$ and $\eta^{i}
=\sigma^{\alpha\beta}_{ik}(\bar t,\bar y)$ implies that
$$
a_{ij}^{\alpha\beta}(\bar t,\bar x)D_{ij}\bar u(\bar t,\bar x)
\leq  a_{ij}^{\alpha\beta}(\bar t,\bar y)
D_{ij}\bar v(\bar t,\bar y)+I+J,
$$
where
$$
I=a_{ij}^{\alpha\beta}(\bar t,\bar x)D_{ij}
\bar u^{(\varepsilon)}(\bar t,\bar x)
-a_{ij}^{\alpha\beta}(\bar t,\bar y)D_{ij}
\bar u^{(\varepsilon)}(\bar t,\bar y),
$$
$$
J=2\nu K\sum_{i,k=1}^{d}|\sigma^{\alpha\beta}_{ik}(\bar t,\bar x)
-\sigma^{\alpha\beta}_{ik}(\bar t,\bar y)|^{2}.
$$
Since $a^{\alpha\beta}$ is uniformly nondegenerate,
 its square root possesses the same smoothness
properties as $a^{\alpha\beta}$ and 
$$
J\leq NK |\bar x-\bar y|^{2\gamma}  .
$$
We now use  \eqref{3.27.6} to get that
$$
 |\bar x-\bar y|^{2\gamma}
\leq  (M/K)^{2\gamma}
N\varepsilon_{0}^{2\gamma\chi}K^{-\chi(1-\gamma)} .
$$
It turns  out that
$$ 
 -2\gamma-\chi(1-\gamma)=-1-\chi/4.
$$
  It follows that (recall that $M\geq 1$)
$$
 |\bar x-\bar y|^{2\gamma}
\leq N M^{2}K^{-1-\chi/4}
 \varepsilon_{0}^{2\gamma\chi} ,
$$
$$
J\leq NM^{2 }K^{-\chi/4} \varepsilon_{0}^{2\gamma\chi} .
$$

  Also note that
$$
I\leq
[a_{ij}^{\alpha\beta}(\bar t,\bar x)-
a_{ij}^{\alpha\beta}(\bar t,\bar y)]
D_{ij}\bar u^{(\varepsilon)}(\bar t,\bar x)
+N
|D^{2}\bar u^{(\varepsilon)}(\bar t,\bar x)-
D^{2}\bar u^{(\varepsilon)}(\bar t,\bar y)|
$$
$$
\leq NM\varepsilon^{\chi-1}
 |\bar x-\bar y|^{ \gamma} 
+NM \varepsilon^{\chi-2} |\bar x-\bar y| =:I_{1}+I_{2},
$$
where the last inequality is obtained by the mean-value theorem
relying on the fact that $|\bar x-\bar y|\leq\varepsilon/2$
 and $\bar x,\bar y
\in G_{2\varepsilon}$,
so that the straight segment connecting these points
lies inside $G_{\varepsilon}$. 
By looking at the estimate of $J$ we get
$$
I_{1}\leq NM^{2}\varepsilon ^{\chi-1} 
  K^{-1/2-\chi/8}.
$$
An easy computation shows that
$$
1/2+\chi/8-(1-\chi)(1-\gamma)/(2\gamma)=\chi/4+\theta_{1},
$$
where $\theta_{1}=(1-\gamma)(8\gamma)^{-1} \chi >0$.
Hence,
$$
I_{1}\leq NM^{2}K^{-\chi/4}[\varepsilon_{0}
^{\chi-1}K^{-\theta_{1}} ].
$$
 
Next,
$$
I_{2}\leq NM \varepsilon_{0}^{2\chi-2}(M/K)K^{(2-2\chi)
(1-\gamma)/
(2\gamma)}.
$$
Here,
$$
(2-2\chi)(1-\gamma)/
(2\gamma)-1=-\chi/4-2\theta_{1},
$$
so that
$$
I_{2}\leq NM^{2}K^{-\chi/4}[
\varepsilon_{0}^{2\chi-2}K^{-2\theta_{1}}].
$$

It turns out that for
\begin{equation}
                                                \label{4.23.6}
K\geq\varepsilon_{0}^{ -\theta_{2}},
\end{equation}
where 
  $\theta_{2}= 
(2\gamma\chi-\chi+1)/\theta_{1}$ we have $\varepsilon_{0}
^{\chi-1}K^{-\theta_{1}},\varepsilon_{0}^{2\chi-2}K^{-2\theta_{1}}
\leq \varepsilon_{0}^{2\gamma\chi}$ and hence
\begin{equation}
                                              \label{4.23.7}
J,I_{1},I_{2}\leq NM^{2}K^{-\chi/4}
\varepsilon_{0}^{2\gamma\chi}.
\end{equation}

Also \eqref{3.31.2} reads
$$
\partial_{t}\bar u (\bar t,\bar x)
\leq\partial_{t}\bar v (\bar s,\bar y)
+\partial_{t}\bar u^{(\varepsilon)}(\bar t,\bar x)
-\partial_{t}
\bar u^{(\varepsilon)}(\bar t,\bar y)
$$
and
as is easy to see
$$
|\partial_{t}\bar u^{(\varepsilon)}(\bar t,\bar x)
-\partial_{t}
\bar u^{(\varepsilon)}(\bar t,\bar y)|\leq  
NI_{2}.
$$
 
It follows that, for $K$ satisfying \eqref{4.23.6} (along with 
\eqref{3.30.3}
and, of course, \eqref{5.27.1}),  we have
$$
\partial_{t}\bar u (\bar t,\bar x)+
a_{ij}^{\alpha\beta}(\bar t,\bar x)
D_{ij}\bar u(\bar t,\bar x)
$$
\begin{equation}
                                           \label{4.23.1}
\leq \partial_{t}\bar v (\bar t,\bar y)
+ a_{ij}^{\alpha\beta}(\bar t,\bar y)
D_{ij}\bar v(\bar t,\bar y)+NM^{2}K^{-\chi/4}
\varepsilon_{0}^{2\gamma\chi}.
\end{equation}  

Next,
$$
D_{i}\bar u(\bar t,\bar x)=2\nu K(\bar x^{i}-\bar y^{i})+
D_{i}\bar u^{(\varepsilon)}(\bar t,\bar x),
$$
$$
 D_{i}\bar v(\bar t,\bar y)=
2\nu K(\bar x^{i}-\bar y^{i})
+D_{i}\bar u^{(\varepsilon)}(\bar t,\bar y),
$$
$$
D_{i}\bar u(\bar t,\bar x)-D_{i}\bar v(\bar t,\bar y)=
D_{i}\bar u^{(\varepsilon)}(\bar t,\bar x)-
D_{i}\bar u^{(\varepsilon)}(\bar t,\bar y),
$$
where
$$
|D \bar u^{(\varepsilon)}(\bar t,\bar x)
-D\bar u^{(\varepsilon)}(\bar t,\bar y)|
\leq NM \varepsilon^{\chi-1}|\bar x-\bar y| =N\varepsilon I_{2}.
$$
This and the rough estimate 
$$
|\bar x-\bar y|\leq N\varepsilon_{0}^{\chi}MK^{-1} 
$$
 following from   \eqref{3.27.6} 
 lead  to   
$$
\hat b^{\alpha\beta}_{i}D_{i}\bar u (\bar t,\bar x)-
\hat b^{\alpha\beta}_{i}D_{i}\bar v (\bar t,\bar y)
\leq
 NM^{2}K^{-\chi/4}\varepsilon_{0}^{2\gamma\chi}+NM
\omega( N\varepsilon_{0}^{\chi}MK^{-1}).
$$

Finally,  
$$
\hat f^{\alpha\beta}(\bar t,\bar x)\leq 
\hat f^{\alpha\beta}(\bar t,\bar y)
+NM
\omega( N\varepsilon_{0}^{\chi}MK^{-1}),
$$
$$
-\hat c^{\alpha\beta}\bar u(\bar t,\bar x)+
 \hat  c^{\alpha\beta}\bar v(\bar t,\bar y)=-
\hat  c^{\alpha\beta}(\bar t,\bar x)
[\bar u(\bar t,\bar x)-\bar v(\bar t,\bar y)]
$$
$$
+\bar v(\bar t,\bar y)
[\hat  c^{\alpha\beta}(\bar t,\bar y)-\hat  c^{\alpha\beta}(\bar t,\bar x)]
\leq-\bar c[\bar u(\bar t,\bar x)-\bar v(\bar t,\bar y)]
+NM
\omega( N\varepsilon_{0}^{\chi}MK^{-1}),
$$
where $\bar c=\inf_{G}\psi^{-1}$ and
the last inequality follows from the fact that
$\hat c^{\alpha\beta}\geq \bar c$ and $\bar u(\bar t,\bar x)-
\bar v(\bar t,\bar y)\geq0$
(see \eqref{3.27.4}).

Therefore we infer from \eqref{3.27.8}, \eqref{3.28.1},
 and the last estimates that
$$
0\leq \supinf_{\alpha\in A\,\,\beta\in B}\big[ 
a_{ij}^{\alpha\beta}D_{ij}\bar v+
\hat b^{\alpha\beta}_{i}D_{i}\bar v-\hat c^{\alpha\beta}
\bar v+f^{\alpha\beta}\big)(\bar t,\bar y)  \big]
$$
$$
  -\bar c[\bar u(\bar t,\bar x)-\bar v(\bar t,\bar y)]
+N_{1} M^{2}\varepsilon_{0}^{2\gamma\chi}K^{-\chi/4}
+N_{1}M\omega(N\varepsilon_{0}^{\chi}MK^{-1}),
$$
$$
\bar u(\bar t,\bar x)-\bar v(\bar t,\bar y)\leq 
N_{1} M^{2}\varepsilon_{0}^{2\gamma\chi}K^{-\chi/4}
+N_{1}M\omega(N_{1}\varepsilon_{0}^{\chi}MK^{-1}).
$$
We can certainly assume that $N_{1}\geq1$.
Then
we take $\mu=2N_{1}$ and take $\kappa \in(0,1)$,
depending only on $\tau$, $\delta$, and $d$,
and $\xi\in(0,\infty)$, depending only on
$\tau$, $K_{0}$, $d$,  $\delta$, and the diameter of $G$, such that
for $\varepsilon_{0}=\xi M^{-1/\kappa}$ and all $M\geq1$
we have 
$$
N_{1} M^{2}\varepsilon_{0}^{2\gamma\chi}\leq 1/2
 ,\quad
N_{1}\varepsilon_{0}^{\chi}M\leq M^{-1/\tau}.
$$
Then 
we arrive at a contradiction with 
\eqref{3.27.4} and, 
since now  \eqref{4.23.6} and \eqref{3.30.3}
are satisfied if $K\geq NM^{1/\eta}$ for  appropriate  
$\eta$ and $N$, the theorem is proved.

 \mysection{Proof of Theorem \protect\ref{theorem 4.1.2}}
                                           \label{section 4.3.3}
 In light of Lemma \ref{lemma 4.1.1}
it suffices to prove \eqref{3.30.01}
only for $K\geq N$, where $N$
depends only  on
  $\tau$,  $d$, $K_{0}$,  $\delta$, and   
$Q$ and
 satisfies \eqref{5.27.1} with $N$
in place of $K$.

Fix a sufficiently small $\varepsilon_{0}>0$ such that
$Q_{\varepsilon_{0}}\ne\emptyset$ and for $\varepsilon
\in(0,\varepsilon_{0}]$ define
$$
\xi_{\varepsilon, K} =\partial_{t}u_{  K}^{(\varepsilon)}
+\max(F[u_{ K}^{(\varepsilon)}],
 P[u_{  K}^{(\varepsilon)}]+ K), 
$$
$$
\xi_{\varepsilon,- K} =\partial_{t}u_{- K}^{(\varepsilon)}
+\min(F[u_{-K}^{(\varepsilon)}],
- P[-u_{-K}^{(\varepsilon)}]+ K), 
$$
 in $Q_{\varepsilon_{0}}$, where we use notation
\eqref{5.3.2}. Since the second-order derivatives  
with respect to $x$ and the first derivatives
with respect to $t$
of $u_{\pm K}$   are bounded in $Q_{\varepsilon_{0}}$,
we have $ \xi_{\varepsilon,\pm K} 
\to0$ as $\varepsilon\downarrow0$
in any $\cL_{p}(Q_{\varepsilon_{0}})$ for any $K$. Furthermore,
$ \xi_{\varepsilon,\pm K} $ are
 continuous. Therefore, there exist smooth functions
$ \zeta_{\varepsilon,K} $ such that
$$
-\varepsilon\leq
  \zeta_{\varepsilon,K}-
\min(\xi_{\varepsilon,K},-\xi_{\varepsilon,-K})\leq 0
$$
in $G_{\varepsilon_{0}}$.

By Theorem 6.4.1 of \cite{Kr_85},
for any subdomain $G'$ of $G_{\varepsilon_{0}}$
of class $C^{3}$,  there exists a unique
$w_{\varepsilon,K}\in C^{1,2}([0,T
-\varepsilon^{2}_{0})\times  G')
\cap C([0,T-\varepsilon^{2}_{0}]\times\bar G')$ satisfying
\begin{equation}
                                             \label{4.28.1}
\partial_{t}w_{\varepsilon,K}+
\sup_{a\in\bS_{\hat\delta},|b|\leq K_{0} }[
a_{ij}D_{ij}w_{\varepsilon,K}+b_{i}D_{i}w_{\varepsilon,K}]=
\zeta_{\varepsilon,K}
\end{equation}
in $[0,T-\varepsilon^{2}_{0})\times  G'$ 
with zero boundary condition
on the parabolic boundary of
$[0,T-\varepsilon^{2}_{0})\times  G'$. Obviously,
$$
\partial_{t}(u_{K}^{(\varepsilon)}-w_{\varepsilon,K})+
\max(F[u_{K}^{(\varepsilon)}-w_{\varepsilon,K}],
P[u_{K}^{(\varepsilon)}-w_{\varepsilon,K}]-K) 
$$
$$
\geq \partial_{t}u_{K}^{(\varepsilon)}
+\max(F[u_{K}^{(\varepsilon)}],
P[u_{K}^{(\varepsilon)}]-K)-\partial_{t}w_{\varepsilon,K}
$$
$$
-
\sup_{a\in\bS_{\hat\delta} ,|b|\leq K_{0}}[
a_{ij}D_{ij}w_{\varepsilon,K}+b_{i}D_{i}w_{\varepsilon,K}]
=\xi_{\varepsilon,K}-\zeta_{\varepsilon,K}\geq0,
$$
$$
\partial_{t}(u_{-K}^{(\varepsilon)}+w_{\varepsilon,K})+
\min(F[u_{-K}^{(\varepsilon)}+w_{\varepsilon,K}],
-P[-u_{-K}^{(\varepsilon)}-w_{\varepsilon,K}]+K) 
$$
$$
\leq \partial_{t}u_{-K}^{(\varepsilon)}+
\min(F[u_{-K}^{(\varepsilon)}],
-P[-u_{-K}^{(\varepsilon)}]+K)+\partial_{t}w_{\varepsilon,K}
$$
$$
+\sup_{a\in\bS_{\hat\delta},|b|\leq K_{0} }[
a_{ij}D_{ij}w_{\varepsilon,K}+b_{i}D_{i}w_{\varepsilon,K}]
=\xi_{\varepsilon,-K}+\zeta_{\varepsilon,K}\leq 0
$$
in $[0,T-\varepsilon^{2}_{0})\times  G')$.
After setting
$$
\mu_{\varepsilon,K}=\sup_{\partial'
([0,T-\varepsilon^{2}_{0})\times  G')}(
u_{K}^{(\varepsilon)}-u_{-K}^{(\varepsilon)}-2w_{\varepsilon,K})_{+}
$$
we conclude by Theorem \ref{theorem 3.27.1} applied to
$u_{K}^{(\varepsilon)}-w_{\varepsilon,K}$
and $u_{-K}^{(\varepsilon)}+w_{\varepsilon,K}+\mu_{\varepsilon,K}$
in place of $u$ and $v$, respectively, that
there exist  a  constant    $N \in(0,\infty)$,
depending only on $\tau$, the diameter of $G$,   
    $d$, $K_{0}$, and $\delta$,
 and a  constant  $\eta\in(0,1)$,
depending only on $\tau$,
    $d$,   and $\delta$,
such that, if $K\geq NM_{\varepsilon,K}^{1/\eta}$, then
$$
u_{K}^{(\varepsilon)}-u_{-K}^{(\varepsilon)}
\leq  \mu_{\varepsilon,K}+w_{\varepsilon,K}+
 NK^{-\chi/4}+  NM\omega(M^{-1/\tau}K^{-1})
$$
 in $[0,T-\varepsilon^{2}_{0})\times  G')$, where 
$M_{\varepsilon,K}\geq1$ is any number satisfying
$$
M_{\varepsilon,K}\geq\|u_{K}^{(\varepsilon)}-w_{\varepsilon,K},
u_{-K}^{(\varepsilon)}+w_{\varepsilon,K}+\mu_{\varepsilon,K}
\|_{C^{1+\chi}([0,T-\varepsilon^{2}_{0})\times  G'))}.
$$

First we discuss
what is happening as $\varepsilon\downarrow0$. By
 the $W^{1,2}_{p}$-theory
(see, for instance, Theorem 1.1 of \cite{DLK_12}) 
$w_{\varepsilon,K}\to 0$ in any $W^{1,2}_{p}$, which by
embedding theorems implies that $w_{\varepsilon,K}\to 0$ in
$C^{1+\chi}([0,T-\varepsilon^{2}_{0})\times  G'))$. Obviously, the constants 
$\mu_{\varepsilon,K}$ converge in 
$C^{1+\chi}([0,T-\varepsilon^{2}_{0})\times  G'))$ to
$$
\sup_{\partial'
([0,T-\varepsilon^{2}_{0})\times  G')}(
u_{K} -u_{-K}  )_{+}.
$$
Now Theorem \ref{theorem 3.31.1} implies that
for sufficiently small $\varepsilon$ one can take
$N\varepsilon_{0}^{-1-\chi}$ as $M_{\varepsilon,K}$,
where $N$ depends only on $d$, $\delta$, $G$, and $K_{0}$.
Thus for sufficiently small $\varepsilon$,
if $K\geq N\varepsilon_{0}^{-(1+\chi)/\eta}$,
then
$$
u_{K}^{(\varepsilon)}-u_{-K}^{(\varepsilon)}
\leq  \mu_{\varepsilon,K}+w_{\varepsilon,K}+
 NK^{-\chi/4}+N\varepsilon_{0}^{-1-\chi}
\omega(\varepsilon_{0}^{  (1+\chi)/\tau}K^{-1})
$$
in $[0,T-\varepsilon^{2}_{0})\times  G'$, which after letting $\varepsilon\downarrow0$
yields
$$
u_{K} -u_{-K} 
\leq   
 NK^{-\chi/4}+N\varepsilon_{0}^{-1-\chi}
\omega(\varepsilon_{0}^{ (1+\chi)/\tau}K^{-1})
$$
in $[0,T-\varepsilon^{2}_{0})\times  G'$.
The arbitrariness of $G'$ and Lemma \ref{lemma 4.1.1}
now allow us to conclude
that for any   $\varepsilon_{0}>0$, for which
$Q_{\varepsilon_{0}}\ne\emptyset$,
$$
u_{K}-u_{-K}\leq  NK^{-\chi/4}
+N\varepsilon_{0}^{-1-\chi}
\omega(\varepsilon_{0}^{ (1+\chi)/\tau}K^{-1})
+\sup_{Q\setminus 
 Q _{\varepsilon_{0}}}(
u_{K} -u_{-K}  )_{+}
$$
\begin{equation}
                                                  \label{4.6.2}
\leq  NK^{-\chi/4}
+N\varepsilon_{0}^{-1-\chi}
\omega(\varepsilon_{0}^{(1+\chi)/\tau}K^{-1})+N\varepsilon_{0}
\end{equation}
in $Q$ provided that  
\begin{equation}
                                                  \label{4.6.3}
K\geq   N_{1}\varepsilon_{0}^{-(1+\chi)/\eta}.
\end{equation}
This obviously proves the 
first assertion of the theorem because as is noted in the proof
of Theorem \ref{theorem 4.1.3} we have $u_{-K}\leq u_{K}$.

To prove the second assertion observe that
  we can certainly assume that
\eqref{5.27.1} holds with $N_{1}$ in place of $K$
 and note that for $\omega=
t^{\tau}$ and $\varepsilon_{0}=K^{-\eta/(1+\chi)}$
condition \eqref{4.6.3} becomes 
$K\geq  N_{1}$ and \eqref{4.6.2}
becomes
$$
u_{K}-u_{-K}\leq  N K^{-\chi/4}+N K^{-\tau} +N
K^{-\eta/(1+\chi)}.
$$
This yields  the desired result and proves the theorem.

\mysection{Proof of assertion (ii)
of Theorem \protect\ref{theorem 4.1.3}}
                                           \label{section 4.3.4}

First we prove uniqueness. Let $w$
be a continuous in $\bar Q$
 viscosity solution of $\partial_{t}w+F[w]=0$ with boundary data
$g$. Observe that in the notation
from Sections \ref{section 4.3.2}
 and \ref{section 4.3.3} we have
$$
\partial_{t}u_{K}^{(\varepsilon)}+
F[u_{K}^{(\varepsilon)}+\bar w_{\varepsilon,K}+\kappa \psi]<0
$$
in $[0,T-\varepsilon^{2})\times G'$ for any $\kappa>0$,
where $\bar w_{\varepsilon,K}$ is a solution of
class $C^{1,2}([0,T
-\varepsilon^{2}_{0})\times  G')
\cap C([0,T-\varepsilon^{2}_{0}]\times\bar G')$
of \eqref{4.28.1} in the domain $[0,T
-\varepsilon^{2}_{0})\times  G'$
with zero   condition on its parabolic boundary
and
 with a smooth $\zeta_{\varepsilon,K}$
satisfying
$$
-\varepsilon\leq \xi_{\varepsilon,K}+
\zeta_{\varepsilon,K}\leq0.
$$

 This and the definition of viscosity solution
imply that the minimum of
$u_{K}^{(\varepsilon)}+\bar w_{\varepsilon,K}+\kappa \psi
-w$ in $[0,T-\varepsilon^{2}]\times \bar G'$ is either positive or is attained
on the parabolic boundary 
of $[0,T-\varepsilon^{2})\times G'$. 
The same conclusion holds after letting
$\varepsilon,\kappa\downarrow0$ and replacing $G'$
with $G_{\varepsilon_{0}}$. Hence, in $Q$
$$
u_{K}-w\geq -\sup_{Q\setminus Q_{\varepsilon_{0}}}
|u_{K}-w|,
$$
which after letting $\varepsilon_{0}\downarrow0$
and then $K\to\infty$
yields $w\leq v$. By comparing $v$ with $u_{-K}$ we get
$w\geq v$, and hence uniqueness.

To prove that $v$ is a viscosity solution we need
the following
Lemma 6.1 of \cite{Kr_14_0} derived there from
Theorem 3.1  of \cite{Kr78} or Theorem 3.3.9
of \cite{Kr_85}.  Introduce
$$
F_{0}(u_{ij},t,x)=F(u_{ij},Dv(t,x),v(t,x),t,x)
$$
$$
C_{r}(t,x)=(0,r^{2})\times\{y\in\bR^{d}:|y|<r\}+(t,x),\quad
C_{r}=C_{r}(0,0).
$$

\begin{lemma}
                              \label{lemma 9.20.1}
There is a constant $N$, depending only on $d$ and
$\delta$, such that
for any $C_{r}(t,x)$ satisfying $C_{r}(t,x)\subset Q$  and
$\phi\in W^{1,2}_{d+1}(C_{r}(t,x))
\cap C(\bar C_{r}(t,x))$ we have on $\bar C_{r}(t,x)$ that
\begin{equation}
                                             \label{9.20.1}
v\leq \phi+Nr^{d/(d+1)}\|(\partial_{t}\phi+
F_{0}[\phi])_{+}\|_{L_{d+1}(C_{r}(t,x))}
+\max_{\partial'C_{r}(t,x)}(v-\phi)_{+} .
\end{equation}
\begin{equation}
                                             \label{9.20.2}
v\geq \phi-Nr^{d/(d+1)}\|(\partial_{t}\phi+
F_{0}[\phi] )_{-}\|_{L_{d+1}(C_{r}(t,x))}
-\max_{\partial'C_{r}(t,x)}(v-\phi)_{-} .
\end{equation}
\end{lemma}

Now let $\phi\in C^{1,2}([0,T]\times\bar G)$ and suppose that 
$v-\phi$ attains a local
maximum at $ (t_{0},x_{0})\in [0,T)\times G$. Without losing generality
we may assume that $t_{0}=0,x_{0}=0$, $w(0)-\phi(0)=0$.
  Then for   $\varepsilon
>0$ and all small $r>0$ for
$$
\phi_{\varepsilon,r}(t,x)=\phi (t,x)+\varepsilon(
|x |^{2}+t - r^{2})
$$
 we have that 
$$
\max_{\partial'C_{r} }(v -\phi_{\varepsilon,r})_{+}
=0.
$$
Hence, by Lemma \ref{lemma 9.20.1}
$$
   \varepsilon r^{2}= 
(v -\phi_{\varepsilon,r})(0)
\leq Nr^{d/(d+1)}\|(\partial_{t}\phi_{\varepsilon,r}+
F_{0}[\phi_{\varepsilon,r}])_{+}\|_{L_{d+1}(C_{r} )}.
$$
It follows that
$$
\sup_{C_{r}}[\partial_{t}\phi_{\varepsilon,r}+
F_{0}[\phi_{\varepsilon,r}]]>0,
$$
which by letting first $r\downarrow0$ and then $\varepsilon
\downarrow0$ yields
$$
0\leq \partial_{t}\phi( 0)+
F_{0}[\phi]( 0)=
\partial_{t}\phi( 0)+
F(D_{ij}\phi( 0),D_{i}\phi( 0),\phi( 0), 0),
$$
where the equality follows from the fact that at $0$
the derivatives of $v-\phi$ with respect to $x$
vanish.
We have  just proved that $v$ is a viscosity subsolution.

Similarly by using \eqref{9.20.2} one
proves that $v$ is a viscosity supersolution.
This proves the theorem.

\mysection{Proof of Theorem \protect\ref{theorem 4.1.4}}
                                          \label{section 5.3.1}

We start with the boundary behavior of $v_{h}$.
\begin{lemma}
                                          \label{lemma 5.1.1}
There exist a constant $N$ such that for all sufficiently
small $h>0$ we have $|v_{h}|\leq N\rho$ in 
$Q_{(h)}$.
\end{lemma}

This lemma is easily proved by using standard barriers
(see, for instance, Lemma 8.8 in \cite{Kr_11}).

We also need the following combination of
Theorems 1.9 and 2.3 of \cite{Kr_12_1}, which provide
a parabolic version of the Fang Hua Lin estimate
and in which 
by $\frL_{\delta,K_{0}}$ we denote the set of parabolic
operators of the form
\begin{equation}
                                                  \label{9.5.7}
L=\partial_{t}+a^{ij}(t,x)D_{ij}+b^{i}(t,x)D_{i}-c(t,x),
\end{equation}
with the coefficients satisfying Assumptions
\ref{assumption 1.9.1} (ii) and (iv) and with $c\geq0$.
Recall that $C_{r}$ are introduced before
Lemma \ref{lemma 9.20.1}.

\begin{theorem}
                       \label{theorem 8.16.1}
Let  
$u\in   C(\bar C_{1})\cap  W^{1,2}_{d+1,loc}(C_{1})$. Then there are
constants $\theta_{0}\in(0, 1]$
 and $N$,
depending only on
$\delta, d$, and $K_{0}$, such that
for any $\theta\in(0,\theta_{0}]$ and 
$L\in  \frL_{\delta,K_{0}} $ we have
\begin{equation}
                                                 \label{8.11.01}
\int_{C_{1}} [|D^2u|^{\theta }+|Du|^{\theta }] \, dx    dt
 \leq  N  \sup_{\partial'C_{1}} 
 |u|^{\theta } + N 
\left(\int_{C_{1}}|  Lu|^{d+1} \,
dx  dt\right)^{\theta/{(d+1)}}.
\end{equation}  
\end{theorem}

In \cite{Kr_12_1} estimate \eqref{8.11.01} is derived
only with $\theta=\theta_{0}$. For $\theta\in(0,\theta_{0}]$
it is obtained by using H\"older's inequality.

\begin{corollary}
                                     \label{corollary 10.7.1}
There exists a constant $\theta_{0}\in(0,1]$,
depending only on $\delta$, $K_{0}$, $d$, and $G$, and
there exists a constant $N$,
depending only on $\delta$, $K_{0}$, $d$, $T$, and $G$, such that
for any  $L\in\frL_{\delta,K_{0}}$,  
 $\theta\in(0,\theta_{0}]$, and $u\in 
W^{ 1, 2}_{d,loc}(Q)\cap C(\bar Q)  $
 we have
\begin{equation}
                                       \label{9.5.06}
 \int_{Q} [|D^2u|^\theta+|Du|^\theta]
\,dxdt \leq N\|Lu\|_{L_{d+1}( Q )}^{\theta}
 +N\sup_{\partial' Q}|u|^{\theta} . 
\end{equation}

\end{corollary}
 Indeed, one can represent $\bar Q$ as the finite union
of the closures of cylinders of height one
in the $t$ variable with bases of class $C^{2}$ each of which
(bases)
admits a one-to-one $C^{2}$ mapping on $\{x:|x|<1\}$
with $C^{2}$ inverse. Then after changing coordinates
one can use Theorem \ref{theorem 8.16.1} applied to appropriately 
changed operator $L$. For the transformed operator the
constants $\delta$ and $K_{0}$
may change but still will only depend
on $\delta, K_{0},d$, and $G$. Then after combining the results
of application of Theorem \ref{theorem 8.16.1} one
obtains \eqref{9.5.06} with $\bar Q$ in place of
$\partial' Q$. However, the parabolic Alexandrov estimate
shows that this replacement can be avoided on the account
of, perhaps, increasing the first $N$ on the right
in \eqref{9.5.06}.
 
The following result is one of our main technical tools.
Everywhere below by $N$ we denote generic constants
independent of $\varepsilon,\varepsilon_{0},K$, $h$,
and the arguments of functions under
consideration.

\begin{theorem}
                                  \label{theorem 1.9.1}
There exist   constants $N\in(0,\infty)$
 and $\theta_{1}\in(0,1)$ such that, for all sufficiently
small $\varepsilon_{0}>0$, for any 
$\varepsilon\in[0,\varepsilon_{0}/2]$, $K\geq1$ 
 and $|l|=1$, we have
\begin{equation}
                                                    \label{5.1.4}
 \|u_{K}(\varepsilon^{2}+\cdot,\varepsilon l+\cdot)
-u_{K}\|_{W^{1,2}_{d+1}
(Q_{ \varepsilon_{0}})}^{d+1} \leq
N\varepsilon^{ \theta_{1} } ( K\varepsilon^{-1}_{0})^{ d+1 }.
\end{equation}
\end{theorem}

To specify what we mean by ``sufficiently
small $\varepsilon_{0}>0$'' let us
  say that a number $\varepsilon >0$
is sufficiently small if Corollary \ref{corollary 10.7.1}
holds with the same $\theta_{0}$,
a constant $N$ which is twice the constant
from \eqref{9.5.06}, and with $Q_{ \varepsilon }$
in place of $Q$. The fact that the set of sufficiently
small $\varepsilon $ contains 
$[0,\alpha)$ with $\alpha>0$ follows from the way
Corollary \ref{corollary 10.7.1} is proved and from the fact
that the boundaries of $G_{\varepsilon}$ have the same
regularity as that of $G$ if $\varepsilon$ is small enough.

Proof. Set $w:=u_{K}(\varepsilon^{2}+\cdot,\varepsilon l+\cdot)
-u_{K}$ and observe that   $\cL_{d+1}(Q_{ \varepsilon_{0}})$-norm 
of $w$ and its sup norm is easily estimated
since $|\partial_{t}u_{K}|+|Du_{K}|\leq NK$
in $Q$. Next, as is easy to see,
there is an operator $L\in\frL_{\hat{\delta},K_{0}}$
such that
\begin{equation}
                                             \label{5.1.6}
Lw+f=0
\end{equation}
in $Q_{ \varepsilon}$, where
$$
f(t,x)=\max(F[u_{K}], P[u_{K}]-K)
( s,y)
$$
$$
-\max(F(D^{2}u_{K}(s,y),D u_{K}(s,y),u_{K}(s,y),t,x),
P[u_{K}]
( s,y)-K),
$$
$s=\varepsilon^{2}+t$, $y=\varepsilon l+x$. It follows that
the estimate of $\partial_{t}w$ can be obtained from
\eqref{5.1.6} once $D^{2}w$, $Dw$, and $f$ are properly estimated.
By the way,  since
$|D^{2}u_{K}|,|Du_{K}|\leq NK\varepsilon^{-1}_{0}$ in 
$Q_{\varepsilon_{0}/2}$ 
and the data $a,b,c,f$ are H\"older continuous we have that
$|f|\leq N\varepsilon^{\tau_{1}}(K\varepsilon^{-1}_{0})$
in $Q_{\varepsilon_{0}}$,
where $\tau_{1}=\min(\tau, 2\gamma_{t},\gamma)$.

We now apply Corollary \ref{corollary 10.7.1}
 to $Q_{ \varepsilon_{0}}$ and $w$
in place of $Q$ and $u$, respectively.
We also use the inequalities like that $|u|^{d+1}
\leq|u|^{\theta_{0}}\sup|u|^{d+1-\theta_{0}}$
while estimating the left-hand side of \eqref{5.1.4}.
Then we obtain 
$$
 \|D^{2}[u_{K}(\varepsilon^{2}+\cdot,\varepsilon l+\cdot)
-u_{K}]\|_{\cL_{d+1}
(Q_{ \varepsilon_{0}})}^{d+1}
$$
$$
+\|D [u_{K}(\varepsilon^{2}+\cdot,\varepsilon l+\cdot)
-u_{K}]\|_{\cL_{d+1}
(Q_{ \varepsilon_{0}})}^{d+1} \leq
N\varepsilon^{ \theta_{1} } ( K\varepsilon^{-1}_{0})^{ d+1 }+I
$$
with $\theta_{1}
=\tau_{1}\theta_{0}$, where $I$ is the boundary term.
However this term is dominated by the right-hand side of
\eqref{5.1.4} due to Theorem \ref{theorem 3.31.1}.
The theorem is proved.

\begin{lemma}
                                   \label{lemma 1.9.2}
There exists a  constant  $N\in(0,\infty)$
  such that, if 
 $\varepsilon_{0}>0$ is sufficiently small,
then for
any $h,\varepsilon\in(0,\varepsilon_{0}/4]$ and $K\geq1$
there exists   $t_{0},x_{0}$ with $|x_{0}|
\leq h$ and $0\leq t_{0}\leq h^{2}$ for which
$$
\sum_{(t,x)\in   Q_{(h)}\cap
Q_{\varepsilon_{0}}}I^{d+1}_{+}(K,\varepsilon,t+t_{0},x+x_{0},t,x)
h^{d+2}
$$
\begin{equation}
                                        \label{1.9.3}
\leq N(h^{ \tau_{1}(d+1)}+\varepsilon^{ \theta_{1} })
(K \varepsilon_{0}^{-1})^{ d+1},
\end{equation}
where  
$$
I (K,\varepsilon,s,y,t,x)
=\partial_{t}u^{(\varepsilon)}_{K}(s,y)
+F(u^{(\varepsilon)}_{K}(s,y),
Du^{(\varepsilon)}_{K}(s,y),D^{2}
u^{(\varepsilon)}_{K}(s,y),t,x) 
$$
and  we use notation
\eqref{5.3.2}.
\end{lemma}

Proof.  Notice that in the left-hand side
of \eqref{1.9.3} similarly to the above proof
$$
|I (K,\varepsilon,t+t_{0},x+x_{0},t,x)-
\partial_{t}u^{(\varepsilon)}_{K}(t+t_{0},x+x_{0})
$$
$$
-
F[u^{(\varepsilon)}_{K}](t+t_{0},x+x_{0})|
\leq Nh^{\tau_{1}}(K\varepsilon^{-1}_{0}).
$$
This contributes a part of the right-hand side of \eqref{1.9.3}.
Concerning the remaining part observe
that $\partial_{t}u_{K}+F[u_{K}]\leq0$ and therefore
$$
\partial_{t}u^{(\varepsilon)}_{K}(t+t_{0},x+x_{0})
+F[u^{(\varepsilon)}_{K}](t+t_{0},x+x_{0})
\leq\partial_{t}[u^{(\varepsilon)}_{K}-u_{K}](t+t_{0},x+x_{0})
$$
$$
+F[u^{(\varepsilon)}_{K}](t+t_{0},x+x_{0})-
F[u _{K}](t+t_{0},x+x_{0})\leq NJ_{\varepsilon,K}(t+t_{0},x+x_{0}),
$$
where
$$
J_{\varepsilon,K}=|\partial_{t}u^{(\varepsilon)} _{K}-
\partial_{t}u_{K}|
+|D^{2}u^{(\varepsilon)} _{K}-
D^{2}u_{K}|+|D u^{(\varepsilon)} _{K}-
D u_{K}| +| u^{(\varepsilon)} _{K}-
 u_{K}|.
$$

Notice that ($C_{r}(t,x)$ are introduced
before Lemma \ref{lemma 9.20.1})
$$
\sum_{(t,x)\in   Q_{(h)}\cap
Q_{\varepsilon_{0}}}I_{C_{h/2}(t,x)}\leq 
I_{Q_{\varepsilon_{0}/2}}
$$
implying that
$$
\dashint_{ C_{h/2}}\sum_{(t,x)\in   Q_{(h)}\cap
Q_{\varepsilon_{0}}}
J^{d+1}_{\varepsilon,K}(t+t_{0},x+x_{0})h^{d+2}
\,dx_{0}dt_{0}
$$
$$
=N\sum_{(t,x)\in   Q_{(h)}\cap
Q_{\varepsilon_{0}}}\int_{ C_{h/2}}
J^{d+1}_{\varepsilon,K}(t+t_{0},x+x_{0}) 
\,dx_{0}dt_{0}
$$
\begin{equation}
                                                \label{5.2.1}
=N\sum_{(t,x)\in   Q_{(h)}\cap
Q_{\varepsilon_{0}}}\int_{ C_{h/2}(t,x)}
J^{d+1}_{\varepsilon,K}( t_{0}, x_{0}) 
\,dx_{0}dt_{0}\leq N\int_{Q_{\varepsilon_{0}/2}}
J^{d+1}_{\varepsilon,K}\,dxdt.
\end{equation}

Furthermore, for $(t,x)\in Q_{\varepsilon_{0}/2}$
we have
$$
|D^{2}u^{(\varepsilon)}(t,x)-D^{2}u(t,x)|^{d+1}\leq
N\dashint_{C_{\varepsilon}}|D^{2}u(t+s,x+y)-D^{2}u(t,x)
|^{d+1}\,dyds.
$$
Similar relations are true for the first
order derivatives and functions themselves.
Therefore, in light of Theorem  \ref{theorem 1.9.1}
$$
\int_{Q_{\varepsilon_{0}/2}}
J^{d+1}_{\varepsilon,K}\,dxdt
\leq N\dashint_{C_{\varepsilon}}
\|
u(\cdot+s,\cdot+y)-u\|_{W^{1,2}_{d+1}(Q_{\varepsilon_{0}/2})}^{d+1}
\,dyds
$$
$$
\leq
N_{1}\varepsilon^{\theta_{1}}
(K\varepsilon_{0}^{-1})^{d+1}.
$$
We conclude  that the average of the first integrand
in \eqref{5.2.1} over $C_{h/2}$ is less than 
$N_{1}\varepsilon^{\theta_{1}}
(K\varepsilon_{0}^{-1})^{d+1}$ which implies that there is
a point $(t_{0},x_{0})\in C_{h/2}$ at which the integrand
itself is less than this quantity and this brings the proof
of the lemma to an end.

{\bf Proof of Theorem  \ref{theorem 4.1.4}}.
Take an $\epsilon  >0$,
set $\lambda=\max\{|l_{k}|\}$, and
observe that for any smooth function $w(t,x)$
given 
in $Q_{ \epsilon }$, we have
in $Q_{ \epsilon  +\lambda h}$  that
$$
\big|\partial_{h,t}w(t,x)+F_{h}[w](t,x)-
\partial_{t}w(t,x)-F[w](t,x)\big|
$$
$$
\leq
N h\sup_{Q_{ \epsilon }}\big[
|\partial_{t}^{2}w|+|D^{3}w|+|D^{2}w|+|Dw|\big],
$$
provided that $Q_{ \epsilon 
+\lambda h}\ne\emptyset$.
We apply  this to $w(t,x)=u_{K}^{(\varepsilon)}(t+t_{0},
x+x_{0})$ for $h,\varepsilon,\varepsilon_{0},t_{0},x_{0}$
from Lemma \ref{lemma 1.9.2} also satisfying
$  \lambda  h\leq\varepsilon_{0}/2$
 and $\epsilon=\varepsilon_{0}-
\lambda h$. Below we only consider $K\geq1$.
Denoting
$$
f(K,\varepsilon,t,x)
=\partial_{h,t}u_{K}^{(\varepsilon)}(t+t_{0},
x+x_{0})+F_{h}[u_{K}^{(\varepsilon)}(\cdot+t_{0},
\cdot+x_{0})](t,x),
$$
 we conclude from Lemma \ref{lemma 1.9.2} that
$$
\sum_{(t,x)\in   Q_{(h)}\cap
Q_{\varepsilon_{0}}}f^{d+1}_{+}(K,\varepsilon,t,x)h^{d+2}
\leq N(h^{ \tau_{1}(d+1)}+\varepsilon^{ \theta_{1} })
(K \varepsilon_{0}^{-1})^{ d+1}
$$
$$
+N h^{d+1} \sup_{Q_{\varepsilon_{0}/2}}\big[
|\partial_{t}^{2}u_{K}^{(\varepsilon)}|+
|D^{3}u_{K}^{(\varepsilon)}|+|D^{2}u_{K}^{(\varepsilon)}|
+|Du_{K}^{(\varepsilon)}|\big]^{d+1},
$$
where the last term admits a rough estimate  (see
Theorems \ref{theorem 3.30.1})
by
$$
Nh^{d+1}\varepsilon^{-2(d+1)}\sup_{Q_{\varepsilon_{0}/4}}[
|\partial_{t}u_{K}|
+|D^{2}u_{K} |+|D u_{K} |
+| u_{K} |\big]^{d+1}
$$
$$
\leq N h^{d+1}\varepsilon^{-2(d+1)}(K\varepsilon_{0}^{-1})^{d+1}.
$$
Also observe that on the discrete parabolic boundary
of $Q_{(h)}\cap
Q_{\varepsilon_{0}}$ we have $|v_{h}|\leq N\varepsilon_{0}$
and $|u_{K}^{(\varepsilon)}(\cdot+t_{0},
\cdot+x_{0})|\leq N\varepsilon_{0}$ owing to Lemmas
\ref{lemma 4.1.1} and \ref{lemma 5.1.1}.
It follows by the discrete maximum principle
of Kuo and Trudinger \cite{KT_93} applied to
$v_{h}- u_{K}^{(\varepsilon)}(\cdot+t_{0},
\cdot+x_{0})$
that in
$Q_{(h)}\cap
Q_{\varepsilon_{0}}$ we have
$$
v_{h}\leq u_{K}^{(\varepsilon)}(\cdot+t_{0},
\cdot+x_{0})+
N(h^{ \tau_{1}}+\varepsilon^{ \theta_{1}/(d+1)})
 K \varepsilon_{0}^{-1} 
+Nh\varepsilon^{-2}K\varepsilon_{0}^{-1}+N\varepsilon_{0}.
$$
 By Theorem  \ref{theorem 3.30.1}
the quantities $|\partial_{t}u_{K}|$ and
$|Du_{K}|$ are bounded by $NK$ in $Q$. Furthermore,
the straight segment connecting any point  
$x\in G_{\varepsilon_{0}}$ with a point of type
$x+\varepsilon y+x_{0}$, where  
 $|y|\leq 1$  and $|x_{0}|\leq
h$ lies in $G_{\varepsilon_{0}/2}\subset G$. It follows
that
$ u_{K}^{(\varepsilon)}(\cdot+t_{0},
\cdot+x_{0})\leq u_{K}+NK(\varepsilon+h)$ on
$G_{\varepsilon_{0}}$.

Combining this with Theorem \ref{theorem 4.1.3}
yields that for   $K\geq 1$ and $h,\varepsilon,\varepsilon_{0}$
as above
$$
v_{h}\leq v+NK^{-\xi}+N(\varepsilon+h)K\varepsilon^{-1}_{0} +
N(h^{ \tau_{1}}+\varepsilon^{ \theta_{1}/(d+1)})
 K \varepsilon_{0}^{-1} 
+Nh\varepsilon^{-2}K\varepsilon_{0}^{-1}+N\varepsilon_{0},
$$
which now holds not only in $Q_{(h)}\cap
Q_{\varepsilon_{0}}$ but also in $Q_{(h)}$ again in light
of Lemmas
\ref{lemma 4.1.1} and \ref{lemma 5.1.1}.

Now, first we take $\varepsilon=h^{1/3}$. Then we obtain
$$
v_{h}\leq v+NK^{-\xi}+Nh^{\theta_{2}}K\varepsilon_{0}^{-1} 
+N\varepsilon_{0},
$$
where $\theta_{2}=\min(\theta_{1}/(3d+3),\tau_{1})$.
Then we take $\varepsilon_{0}=h^{\theta_{2}/2}$
(and only concentrate on $h$ such that 
$h^{1/3},
  h\leq\varepsilon_{0}/4,
 \lambda h\leq\varepsilon_{0}/2$) and get
$$
v_{h}\leq v+NK^{-\xi}+NKh^{\theta_{2}/2} ,
$$
which for $K=h^{- \eta}$, where $\eta= \theta_{2}/(2+2\xi)$,
finally leads to $v_{h}\leq v+Nh^{\xi \eta}$ in $Q_{h}$.

The reader understands that one can prove the inequality
$v_{h}\geq v-Nh^{\xi \eta}$ in $Q_{h}$ by using $u_{-K}$
in place of $u_{K}$. The theorem is proved.

\end{document}